\documentclass{amsart}
\usepackage{amssymb}

\newtheorem{theorem}{Theorem}
\newtheorem{lemma}[theorem]{Lemma}

\theoremstyle{remark}
\newtheorem*{remark}{Remark}
\newtheorem*{notation}{Notation}

\newcommand{\BC}{\mathbb{C}}
\newcommand{\BF}{\mathbb{F}}
\newcommand{\BQ}{\mathbb{Q}}
\newcommand{\BZ}{\mathbb{Z}}

\newcommand{\calA}{\mathcal{A}}
\newcommand{\calB}{\mathcal{B}}
\newcommand{\calO}{\mathcal{O}}

\newcommand{\frakA}{\mathfrak{A}}
\newcommand{\frakp}{\mathfrak{p}}

\newcommand{\xbar}{\overline{x}}

\newcommand{\pib}{\overline{\pi}}

\renewcommand{\hat}{\widehat}
\newcommand{\fq}{{\BF_q}}

\DeclareMathOperator{\Aut}{Aut}
\DeclareMathOperator{\Cl}{Cl}
\DeclareMathOperator{\End}{End}
\DeclareMathOperator{\Gal}{Gal}
\DeclareMathOperator{\res}{Res}
\DeclareMathOperator{\Spec}{Spec}

\begin{document}

\title[Principally polarizable isogeny classes]
      {Principally polarizable isogeny classes\\
       of abelian surfaces over finite fields}

\author[Howe]{Everett W.~Howe}  
% Don't hbox the first author; it doesn't work here as it does below.
\address{Center for Communications Research,
         4320 Westerra Court,
         San Diego, CA 92121-1967, USA.}
\email{however@alumni.caltech.edu}
\urladdr{http://www.alumni.caltech.edu/\~{}however/}

\author[Maisner]{\hbox{Daniel Maisner}}   
% The hbox inhibits linebreaks.
\address{Universidad Aut{\'o}noma de la Ciudad de M{\'e}xico (UACM),
         Plantel San Lorenzo Tezonco, 
         Colegio de Ciencia y Tecnolog{\'i}a,
         Prolongaci{\'o}n San Isidro \# 151, San Lorenzo Tezonco,
         CP 09790 Iztapalapa, D.F. M{\'e}xico.}
\email{danielm@mat.uab.es}

\author[Nart]{\hbox{Enric Nart}}
\address{Departament de Matem{\`a}tiques,
         Universitat Aut{\`o}noma de Barcelona,
         Edifici C, 08193 Bellaterra, Barcelona, Spain.}
\email{nart@mat.uab.es}

\author[Ritzenthaler]{\hbox{Christophe Ritzenthaler}}
\address{Institut de Math{\'e}matiques de Luminy,
         UMR 6206 du CNRS,
         Luminy, Case 907, 13288 Marseille, France.}
\email{ritzenth@iml.univ-mrs.fr}

\date{27 February 2006}



%%%%%%%%%%%%%%%%%%%%%%%%%%%%%%%%%%%%%%%%%%%%%%%%%%%%%%%%%%%%%%%%%%%%%%%%%%%%%
%%                                                                         %%
%% (2000 Mathematics Subject Classification)                               %%
%%                                                                         %%
%% 11            Number theory                                             %%
%%    11G           Arithmetic algebraic geometry (Diophantine geometry)   %%
%%        11G10        Abelian varieties of dimension > 1                  %%
%%        11G25        Varieties over finite and local fields              %%
%% 14            Algebraic Geometry                                        %%
%%    14G           Arithmetic problems. Diophantine geometry              %%
%%        14G15        Finite ground fields                                %%
%%                                                                         %%
%%%%%%%%%%%%%%%%%%%%%%%%%%%%%%%%%%%%%%%%%%%%%%%%%%%%%%%%%%%%%%%%%%%%%%%%%%%%%

\begin{abstract}
Let $\calA$ be an isogeny class of abelian surfaces over $\fq$ with Weil 
polynomial $x^4 + ax^3 + bx^2 + aqx + q^2$.  We show that $\calA$ does not 
contain a surface that has a principal polarization if and only if 
$a^2 - b = q$ and $b < 0$ and all prime divisors of $b$ are congruent to 
$1$ modulo~$3$.
\end{abstract}

\maketitle

\section{Introduction}
An isogeny class of abelian varieties over a field $k$ is 
\emph{principally polarizable} if it contains a variety that admits a 
principal polarization defined over~$k$.  If $k$ is algebraically closed then 
every isogeny class of abelian varieties over $k$ is principally polarizable, 
but over non-algebraically closed fields there may well exist isogeny classes 
that contain no principally polarized varieties.  In general, it can be hard to 
tell whether an isogeny class of abelian varieties is principally polarizable.  
In this paper we consider the case of abelian surfaces over finite fields,
and we give a simple criterion for deciding whether or not an isogeny class 
of such surfaces is principally polarizable.  We use this result in a 
forthcoming paper~\cite{HNR}, in which we determine which isogeny classes of 
abelian surfaces contain Jacobians.

The \emph{Weil polynomial} of an abelian variety over a finite field is the
characteristic polynomial of its Frobenius endomorphism.  The Honda-Tate 
theorem~\cite{tate} shows that two abelian varieties over a finite field are
isogenous to one another if and only if they share the same Weil polynomial.
If $A$ is an abelian surface over $\fq$ then its Weil polynomial has the form
$$x^4+ax^3+bx^2+aqx+q^2$$
for some integers $a$ and $b$; suppressing the dependence on $q$, we will 
let $\calA_{(a,b)}$ denote the isogeny class of the surface $A$.

\begin{theorem}
\label{T:Main} 
Let $\calA=\calA_{(a,b)}$ be an isogeny class of abelian surfaces over~$\fq$.
Then $\calA$ is not principally polarizable if and only if the following three
conditions are satisfied{\/\rm:}
\begin{itemize}
\item[{\rm (a)}]  $a^2-b=q$,
\item[{\rm (b)}]  $b<0$, and 
\item[{\rm (c)}]  all prime divisors of $b$ are congruent to $1$ modulo $3$.
\end{itemize}
\end{theorem}

Recall that an abelian surface over a field of characteristic~$p$ is
\emph{ordinary} if the rank of its group of geometric $p$-torsion points 
is~$2$; it is \emph{supersingular} if this rank is $0$.  If an abelian surface 
over a field of positive characteristic is neither ordinary nor supersingular,
we will call it \emph{mixed}.

Several special cases of Theorem~\ref{T:Main} appear in the literature or
are easily derivable from known results.  For simple ordinary abelian surfaces,
Theorem~\ref{T:Main} occurs as~\cite[Thm.~1.3]{TAMS}.  A straightforward
argument using~\cite[Lem.~2.4]{mn} shows that the theorem holds
for split isogeny classes.  The second and third authors, using techniques 
from~\cite{JAG}, proved that every mixed isogeny class of abelian surfaces is
principally polarizable~\cite[Thm.~4.3]{mn}, and since an isogeny class 
$\calA_{(a,b)}$ is mixed if and only if $p\nmid a$ and $p\mid b$ (so that 
$a^2 - b\neq q$), Theorem~\ref{T:Main} holds for these isogeny classes as 
well.  For simple supersingular isogeny classes, Theorem~\ref{T:Main} is 
equivalent to the following result.

\begin{theorem}
\label{T:Supersingular} 
Let $\calA$ be an isogeny class of simple supersingular abelian surfaces 
defined over $\fq$, where $q$ is a power of a prime $p$.  Then $\calA$ is not
principally polarizable if and only if $\calA=\calA_{(0,-q)}$ and 
$p\equiv 1\bmod3$.
\end{theorem}

To see that the two theorems are equivalent, we consult the complete list of 
simple supersingular isogeny classes given in~\cite[Thm.~2.9]{mn}, and we 
note that the only class that meets the three conditions of
Theorem~\ref{T:Main} is the isogeny class mentioned in 
Theorem~\ref{T:Supersingular}.

Correctly interpreted, Theorem~\ref{T:Main} reveals a geometric feature of the
isogeny classes that are not principally polarizable.  Let $E$ be an elliptic
curve over a finite field~$k$ and let $\ell$ be the cubic extension of $k$.  
There is a natural trace map from the restriction of scalars $\res(\ell/k,E)$
to $E$; the kernel $A$ of this map is the \emph{trace-zero surface} associated 
to $E$.  (The trace-zero surface is called the \emph{reduced restriction of 
scalars} in~\cite{Texel}.)  If the Weil polynomial for $E$ is $x^2 - ax + q$
then $A$ lies in the isogeny class $\calA_{(a,a^2-q)}$, so the isogeny class 
of $A$ satisfies the first condition of Theorem~\ref{T:Main}.   On the other
hand, it is not hard to show that if an isogeny class $\calA$ satisfies all
three conditions of Theorem~\ref{T:Main}, then $\calA$ contains the trace-zero
surface of an elliptic curve over~$\fq$.  (Indeed, for ordinary isogeny classes
the first of the three conditions is sufficient.)  These observations allow
us to restate Theorem~\ref{T:Main} as follows:

\begin{theorem}
\label{T:TraceZero}
Let $E$ be an elliptic curve defined over $\fq$ and let $a$ be the trace of
Frobenius of $E$.  Suppose that $|a|<\sqrt{q}$ and that all of the prime 
divisors of $a^2-q$ are congruent to $1$ modulo~$3$. Then the isogeny class 
of the trace-zero surface constructed from $E$ is not principally polarizable.
Conversely, every isogeny class of abelian surfaces that is not principally 
polarizable contains the trace-zero surface of an elliptic curve whose trace
satisfies the properties listed above.
\end{theorem}

In particular, we see that isogeny classes of abelian surfaces over finite
fields that are not principally polarizable split over the cubic extension of
the base field.  The analogous statement for isogeny classes of higher
dimension is false; for instance, there are absolutely simple isogeny classes 
of four-dimensional abelian varieties over $\BF_{41}$ that are not principally
polarizable~\cite[Ex.~13.4]{TAMS}.

The techniques we use to prove Theorem~\ref{T:Main} are based on ideas from a 
series of papers~\cite{TAMS, JAG, Texel} in which the first author investigated
the obstruction to the existence of a principally polarized variety in a given
isogeny class over an arbitrary field.  The general machinery provided
in~\cite{JAG} can in principle be used to determine whether or not an isogeny
class of abelian varieties over a finite field is principally polarizable.  In
practice, however, there are certain situations in which one can apply the
machinery only if one has sufficient information about some polarization of 
some variety in the isogeny class.  One such situation arises in this paper.

In Section~\ref{S:Polarizable} we use results from~\cite{JAG} to show that the
simple supersingular isogeny classes not listed in 
Theorem~\ref{T:Supersingular} are principally polarizable.  
In Section~\ref{S:Artin} we prove Theorem~\ref{T:Artin}, which specializes the
results of~\cite{JAG} to the case of certain simple isogeny classes and shows
how the value of the Artin symbol of an ideal obtained from a polarization 
determines whether or not the variety is isogenous to a principally polarized
variety.  Finally, in Section~\ref{S:Polarization} we use Theorem~\ref{T:Artin}
to complete the proof of Theorem~\ref{T:Supersingular}  by producing an 
explicit polarization of a trace-zero surface obtained from an elliptic curve 
and computing the Artin symbol of the associated ideal. 

\begin{notation}
Let $k$ be a field and let $\ell$ be a finite extension of $k$.
If $A$ is a variety over $k$ we let $A_\ell$ denote the variety
$A\times_{\Spec k}{\Spec \ell}$ over~$\ell$.  If $A$ is a variety 
over~$\ell$, we let $\res(\ell/k, A)$ denote the Weil restriction of
scalars of $A$ from $\ell$ to~$k$.  Note that if $A$ is an abelian variety 
over~$\ell$ that has a principal polarization, then $\res(\ell/k, A)$ is an
abelian variety over $k$ that has a principal polarization.
\end{notation}

%%%%%%%%%%%%%%%%%%%%%%%%%%%%%%%%%%%%%%%%%%%%%%%%%%%%%%%%%%%%%%%%%%%%%%%%
%%%%%%%%%%%%%%%%%%%%%%%%%%%%%%%%%%%%%%%%%%%%%%%%%%%%%%%%%%%%%%%%%%%%%%%%
\section{Principally polarizable supersingular abelian surfaces}
\label{S:Polarizable}

In this section we prove the `only if' direction of 
Theorem~\ref{T:Supersingular}:  We show that if $\calA$ is an isogeny class of
simple supersingular abelian surfaces over $\fq$ of characteristic~$p$, and if
either $p\not\equiv1\bmod 3$ or $\calA \neq \calA_{(0,-q)}$, then $\calA$ is
principally polarizable.   To prove this, we use one of the main theorems 
of~\cite{JAG} together with a simple result.  We start by proving the simple 
result.

\begin{lemma}
\label{L:Restriction}
Let $A$ be a simple abelian surface defined over $\BF_q$. Then $A_{\BF_{q^2}}$ 
is not simple if and only if $A$ is isogenous to $\res(\BF_{q^2}/\BF_q,E)$ for 
some elliptic curve $E$ defined over~$\BF_{q^2}$.
\end{lemma}

\begin{proof}
Suppose that $A\in\calA_{(a,b)}$ and $A_{\BF_{q^2}}$ is not simple. 
By~\cite[Thm.2.15]{mn} all simple mixed surfaces are absolutely simple; thus, 
$A$ is either ordinary or supersingular. By Lemma~2.1, Proposition~2.14, and 
Table~1 of~\cite{mn}, in both cases we have
\begin{itemize}
\item $a=0$, and
\item $b$ is the trace of Frobenius of an elliptic curve $E$ over $\BF_{q^2}$.
\end{itemize}
It follows from~\cite[\S1(a)]{milne} that the abelian surface 
$\res(\BF_{q^2}/\BF_q, E)$ has the same Weil polynomial as does $A$, so the two
surfaces are isogenous.

Conversely, if $A$ is isogenous to  $\res(\BF_{q^2}/\BF_q,E)$ then 
$A_{\BF_{q^2}}$ is isogenous to the product of $E$ with its $\fq$-conjugate,
so $A$ is not simple over~$\BF_{q^2}$.
\end{proof}

\begin{remark}
The Weil restriction of an elliptic curve is a principally polarizable
variety, so the lemma shows that if $A_{\BF_{q^2}}$ is not simple then
$A$ is isogenous to a principally polarized variety.
\end{remark}

Let us recall the first main theorem from~\cite{JAG}.  Suppose $\calA$ is an
isogeny class of simple abelian varieties over a finite field and let $f$ be
its Weil polynomial.  Since $\calA$ is simple, the polynomial $f$ is a
power of an irreducible polynomial.  Let $\pi$ be a root of $f$ in $\BC$.  The 
field $K = \BQ(\pi)$ is either totally real or a totally imaginary quadratic
extension of a totally real field $K^+$.  In either case, complex conjugation
$x\mapsto \xbar$ induces an automorphism of $K$, trivial if $K$ is totally
real.  

\begin{theorem}[{\cite[Thm.~1.1]{JAG}}]
\label{T:JAG}
If $K$ is totally real then $\calA$ contains a principally polarized variety.
Suppose $K$ is a CM-field.  If a finite prime of $K^+$ ramifies in $K/K^+$, or 
if there is a prime of $K^+$ that divides $\pi-\pib$ and that is inert 
in~$K/K^+$, then $\calA$ contains a principally polarized variety.
\end{theorem}

Now we have the machinery we need to prove the `only if' direction of
Theorem~\ref{T:Supersingular}.  Let $\calA$ be an isogeny class of simple 
supersingular abelian surfaces over~$\fq$.  If the surfaces in $\calA$ split
over $\BF_{q^2}$ then Lemma~\ref{L:Restriction} shows that there is a 
principally polarized surface in $\calA$, so we need only concern ourselves
with the isogeny classes that remain simple over $\BF_{q^2}$.  If the Weil
polynomial of $\calA$ is not irreducible, then the field $K$ has degree at 
most $2$ over $\BQ$, so either $K$ is totally real or $K$ is an imaginary
quadratic extension of $\BQ$ (and hence ramified at a finite prime).  In either
case, Theorem~\ref{T:JAG} shows that $\calA$ is principally polarizable.  Thus 
we need only consider isogeny classes whose Weil polynomials are irreducible. 

Theorem~2.9 of~\cite{mn}, under the heading `(SS1)', lists the irreducible 
supersingular Weil polynomials coming from abelian surfaces.  Table~1 
of~\cite{mn} lists the degrees of the smallest extension fields over which the
isogeny classes split.  Combining these sources, we can make a list of all of 
the Weil polynomials we must consider.  We present this list in 
Table~\ref{Table:WP}.  For all but two of the entries, we give a reason why
there is a principally polarized surface in the corresponding isogeny class.

\begin{table}
\begin{center}
\renewcommand{\arraystretch}{1.1}
\begin{tabular}{|l|l|l|}
\hline
$(a,b)$             & Conditions on $p$, $q$            & Reason for principal polarization \\
\hline\hline
$(0,-q)$            & $p\equiv 1\bmod3$, $q$ nonsquare  & ---\\
\hline
$(0,-q)$            & $p\equiv 7\bmod12$, $q$ square    & ---\\
\hline
$(0,0)$             & $p\equiv 1\bmod4$, $q$ nonsquare  & $K/K^+$ ramified over $2$\\
\hline
$(0,0)$             & $p\equiv 5\bmod8$, $q$ square     & $K/K^+$ ramified over $2$\\
\hline
$(0,q)$             & $p\equiv 1\bmod3$, $q$ nonsquare  & $K/K^+$ ramified over $3$\\
\hline
$(\pm\sqrt q,q)$    & $p\not\equiv 1\bmod5$, $q$ square & $K/K^+$ ramified over $5$\\
\hline
$(\pm\sqrt{5q},3q)$ & $p=5$, $q$ nonsquare              & $K/K^+$ ramified over $5$\\
\hline
$(\pm\sqrt{2q}, q)$ & $p=2$, $q$ nonsquare              & $\pi-\pib$ divisible by inert prime over $2$\\
\hline
\end{tabular}
\end{center} 
\vspace{1ex}
\caption{Certain irreducible supersingular Weil polynomials 
$x^4 + ax^3 + bx^2 + aqx + q^2$ over finite fields $\fq$ of 
characteristic~$p$, together with reasons (if any) why the corresponding
isogeny class of abelian surfaces contains a principally polarized variety.
}
\label{Table:WP} 
\end{table}

We have shown that every isogeny class of simple supersingular abelian surfaces
over $\fq$ other than the first two entries of Table~\ref{Table:WP} is
principally polarizable, which proves the `only if' part of
Theorem~\ref{T:Supersingular}. \qed

\begin{remark}
At first glance it seems as if we have overlooked the case in which $q$ is
a square, $p\equiv 1\bmod 12$, and $\calA = \calA_{(0,-q)}$. But as we see 
from~\cite[Thm.2.9]{mn}, there is no isogeny class with Weil polynomial 
$t^4-qt^2+q^2$ when $q$ is a square and $p\equiv 1\bmod 12$.
\end{remark}

%%%%%%%%%%%%%%%%%%%%%%%%%%%%%%%%%%%%%%%%%%%%%%%%%%%%%%%%%%%%%%%%%%%%%%%%
%%%%%%%%%%%%%%%%%%%%%%%%%%%%%%%%%%%%%%%%%%%%%%%%%%%%%%%%%%%%%%%%%%%%%%%%
\section{Determining the existence of a principal polarization}
\label{S:Artin}

To complete the proof of Theorem~\ref{T:Supersingular}, we must show that when 
$p\equiv 1\bmod 3$ there is no principally polarized abelian surface with Weil
polynomial $x^4 - qx^2 + q^2$.  To accomplish this, we first provide a 
criterion for determining whether or not an isogeny class is principally 
polarizable; the criterion requires sufficient knowledge of an arbitrary 
polarization of a variety in the isogeny class, and is implicit in~\cite{JAG}.  
Then, in Section~\ref{S:Polarization}, we construct a polarization of a 
variety in the isogeny class with Weil polynomial $x^4 - qx^2 + q^2$ to which 
we can apply the criterion.

Our criterion applies to isogeny classes of simple abelian varieties over~$\fq$
of arbitrary dimension.  Let $\calA$ be such an isogeny class, let $A$
be any variety in $\calA$, and let $F$ and $V$ be the Frobenius and 
Verschiebung endomorphisms of $A$.  Then the subring $R:=\BZ[F,V]$ of $\End(A)$
is a domain, and up to isomorphism it is independent of the choice of 
$A\in\calA$.  Let $K := \BQ(F)$ be the quotient field of $R$.  Then $K$ is 
either a totally real field or a CM-field, that is, a totally imaginary 
quadratic extension of a totally real field $K^+$.

As we noted above in Theorem~\ref{T:JAG}, if $K$ is totally real then $\calA$
is principally polarizable; if $K$ is a CM-field that is ramified over $K^+$ at
a finite prime, then $\calA$ is principally polarizable; and if $F-V$ is 
divisible by a prime of $K^+$ that is inert in $K/K^+$, then $\calA$ is
principally polarizable.  Thus, we may henceforth assume that $K$ is a 
CM-field that is unramified over $K^+$ at every finite prime, and that every
prime of $K^+$ that divides $F-V$ is split in $K/K^+$.

Let $\calO$ be the ring of integers of $K$ and let $\calO^+$ be the ring of 
integers of $K^+$.  Let $A$ be any variety in $\calA$ such that 
$\calO\subseteq\End(A)$;  we know from~\cite[Thm.~3.13]{waterhouse} that
such varieties exist.  Let $\hat{A}$ be the dual variety of $A$, and 
suppose there is a polarization $\lambda\colon A\to\hat{A}$ of $A$ whose
degree is coprime to the characteristic of~$k$.  Then the group of 
geometric points of $\ker(\lambda)$ is naturally an $\calO$-module, and the
proof of~\cite[Prop.~7.1]{JAG} shows that there is an ideal $\frakA$ of 
$\calO^+$ such that 
$$\ker(\lambda) \cong \calO/\frakA\calO$$
as $\calO$-modules.  Let $\psi$ denote the Artin map from the ideal group of 
$\calO^+$ to the Galois group of $K$ over~$K^+$; we identify this Galois group
with $\{\pm1\}$.

\begin{theorem}
\label{T:Artin}
The isogeny class $\calA$ is principally polarizable if and only if 
$\psi(\frakA) = 1$.
\end{theorem}

\begin{proof}
In~\cite[\S6]{JAG} the first author defined a contravariant functor $\calB$
from a certain category of rings (the ``proper real/CM-orders'') to the
category of finite $2$-torsion groups.  In particular, for every order $S$ in 
$K$ that is fixed by complex conjugation, we get a group $\calB(S)$, and the 
inclusion map $i\colon R\to \calO$ gives us a group homomorphism 
$$i^*\colon \calB(\calO)\to\calB(R).$$
Furthermore,~\cite[Prop.~6.2]{JAG} shows that 
$$\calB(\calO) \cong \Cl^+(\calO^+)/N_{\calO/\calO^+}(\Cl(\calO)),$$
so that the Artin map gives an isomorphism from $\calB(\calO)$ 
to~$\Gal(K/K^+)$.

The third main result of~\cite{JAG} (Theorem~1.3) says that there is a 
naturally-defined  element $I_\calA$ of $\calB(R)$ that is zero precisely 
when $\calA$ is principally polarizable.  Proposition~7.1 of~\cite{JAG} shows
that $I_\calA$ lies in the image of $\calB(\calO)$ under~$i^*$. In particular,
the proof of~\cite[Prop.~7.1]{JAG} shows that $I_\calA$ is the image under
$i^*$ of the class of the ideal~$\frakA$.  Proposition~7.2 of~\cite{JAG} says
that under our assumptions on $A$ (namely, that $K$ is a CM-field that is 
unramified over $K^+$ at all finite primes, and that no inert prime of $K^+$
divides $F-V$),  the map $i^*$ is an injection.  Therefore, $I_\calA$ is zero
if and only if the class of $\frakA$ in $\calB(\calO)$ is zero, and this 
happens if and only if the Artin symbol on $\frakA$ is trivial.
\end{proof}

\begin{remark}
It is only for convenience that we assume that the polarization $\lambda$ 
of $A$ has degree coprime to the characteristic.  We can associate an 
ideal $\frakA$ of $\calO$ to any polarization of $A$ as in the proof
of~\cite[Prop.~7.1]{JAG}, and Theorem~\ref{T:Artin} remains true
for the ideals associated to these more general polarizations.
\end{remark}

%%%%%%%%%%%%%%%%%%%%%%%%%%%%%%%%%%%%%%%%%%%%%%%%%%%%%%%%%%%%%%%%%%%%%%%%
%%%%%%%%%%%%%%%%%%%%%%%%%%%%%%%%%%%%%%%%%%%%%%%%%%%%%%%%%%%%%%%%%%%%%%%%
\section{An explicit polarization}
\label{S:Polarization}

To use Theorem~\ref{T:Artin} to complete the proof of 
Theorem~\ref{T:Supersingular}, we must construct an explicit polarization of an
abelian surface in the isogeny class with Weil polynomial $x^4 - qx^2 + q^2$.
In this section we construct such a polarization by using ideas
from~\cite{Texel}.  We assume throughout this section that $p\equiv1\bmod3$
and that $p\equiv3\bmod4$ if $q$ is a square.

Our assumptions on $p$ and $q$ guarantee that there is a supersingular 
elliptic curve $E$ over $\fq$ with Weil polynomial $x^2 + q$ 
(see~\cite[Thm.~4.1]{waterhouse}).  Furthermore, 
by~\cite[Thm.~3.13]{waterhouse} we can choose $E$ so that its endomorphism
ring is the maximal order $S$ of the field~$L := \BQ(\sqrt{-q})$.

Let $A$ be the trace-zero surface over $k$ associated to $E$, by which we mean 
the kernel of the trace map from $\res(\BF_{q^3}/\BF_q,E_{\BF_{q^3}})$ to~$E$. 
By~\cite[\S1(a)]{milne}, the Weil polynomial of $A$ is $f = x^4 - qx^2 + q^2$, 
so $A$ lies in the isogeny class $\calA=\calA_{(0,-q)}$ that we are interested
in. Let $F$ and $V$  be the Frobenius and Verschiebung endomorphisms of $A$; 
since $f$ is irreducible, the subring $K = \BQ(F)$ of $\End(A)\otimes\BQ$ is
a field.  In fact, $K$ is a CM-field whose maximal real subfield $K^+$ is
isomorphic to~$\BQ(\sqrt{3q})$.  There are also two imaginary quadratic 
subfields of~$K$; one of them is isomorphic to the endomorphism algebra 
$L = \BQ(\sqrt{-q})$ of $E$, and one is $M:=\BQ(\omega)$, where $\omega$
is a primitive cube root of unity.

The discriminants of $L$ and $M$ are coprime to each other, so the field $K$ 
is unramified over $K^+$ at all finite primes.  We compute that $(F-V)^2 = -q$,
that there is a unique prime $\frakp$ of $K^+$ over~$p$, and that $\frakp$
splits in~$K$.  Thus the field $K$ satisfies all of the hypotheses set forth
at the beginning of Section~\ref{S:Artin}.

In Section~2.2 of~\cite{Texel} it is shown that $A$ is the 
$\BF_{q^3}/\BF_q$-twist of $E\times E$ corresponding to the element of 
$H^1(\Gal(\BF_{q^3}/\BF_q),\Aut(E\times E))$ represented by the cocyle that\
sends the Frobenius automorphism of $\BF_{q^3}/\BF_q$ to the automorphism of
$E\times E$ given by the matrix
$$\zeta = \left[\begin{matrix}-1&-1\\1&0\end{matrix}\right].$$
By replacing $\omega$ with its complex conjugate, if necessary, we may assume 
that $\zeta = \omega$.  Note that every geometric endomorphism of $E\times E$ 
is defined over $\fq$, and that the endomorphism ring of $E\times E$ is 
isomorphic to the ring $M_2(S)$ of $2\times2$ matrices over~$S$. 
Proposition~2.1 of~\cite{Texel} says that an element of $M_2(S)$
gives rise to an endomorphism of $A$ if and only if it commutes with the 
element $\zeta$.  This shows that $\End A \cong S[\omega]$.  Since the
discriminants of $L$ and $M$ are coprime to one another, the ring $S[\omega]$ 
is the full ring of integers $\calO$ of $K$.  Thus we can apply 
Theorem~\ref{T:Artin} to the variety~$A$ if we can find an explicit
polarization of~$A$.

Let $b$ be the endomorphism 
$$\left[\begin{matrix}2&1\\1&2\end{matrix}\right]$$
of $E\times E$, let $\mu$ be the product polarization of $E\times E$, and
let $\lambda = \mu b$, so that $\lambda$ is a polarization of degree~$9$.
As in the proof of~\cite[Lem.~2.6]{Texel}, Proposition~2.2 of~\cite{Texel}
shows that $\lambda$ gives rise to a polarization on $A$.  As we noted at
the beginning of Section~\ref{S:Artin}, there is an ideal $\frakA$ of the ring
of integers $\calO^+$ of $K^+$ such that $\ker(\lambda)\cong\calO/\frakA\calO$,
and the only possibility is that $\frakA$ is the unique prime of $\calO^+$
lying over the rational prime~$3$.

The prime $3$ ramifies in $M/\BQ$ and is inert in $L/\BQ$, so the prime 
$\frakA$ is inert in $K/K^+$.  Therefore its Artin symbol is nontrivial, and 
it follows from Theorem~\ref{T:Artin} that the isogeny class $\calA$ is not
principally polarizable.  This completes the proof of 
Theorem~\ref{T:Supersingular}.\qed

\end{document}